\documentstyle[12pt]{article} 
\begin{document}
\newcommand{\singlespace}{
    \renewcommand{\baselinestretch}{1}
\large\normalsize}
\newcommand{\doublespace}{
   \renewcommand{\baselinestretch}{1.2}
   \large\normalsize}
\renewcommand{\theequation}{\thesection.\arabic{equation}}

\input amssym.def
\input amssym
\setcounter{equation}{0}
\def \ten#1{_{{}_{\scriptstyle#1}}}
\def \Z{\Bbb Z}
\def \C{\Bbb C}
\def \R{\Bbb R}
\def \Q{\Bbb Q}
\def \N{\Bbb N}
\def \l{\lambda}
\def \V{V^{\natural}}
\def \wt{{\rm wt}}
\def \tr{{\rm tr}}
\def \Res{{\rm Res}}
\def \End{{\rm End}}
\def \Aut{{\rm Aut}}
\def \mod{{\rm mod}}
\def \Hom{{\rm Hom}}
\def \im{{\rm im}}
\def \<{\langle} 
\def \>{\rangle} 
\def \w{\omega}
\def \c{{\tilde{c}}}
\def \o{\omega}
\def \t{\tau }
\def \ch{{\rm ch}}
\def \a{\alpha }
\def \b{\beta}
\def \e{\epsilon }
\def \la{\lambda }
\def \om{\omega }
\def \O{\Omega}
\def \qed{\mbox{ $\square$}}
\def \pf{\noindent {\bf Proof: \,}}
\def \voa{vertex operator algebra\ }
\def \voas{vertex operator algebras\ }
\def \p{\partial}
\def \1{{\bf 1}}
\def \ll{{\tilde{\lambda}}}
\def \H{{\bf H}}
\def \F{{\bf F}}
\def \h{{\frak h}}
\def \g{{\frak g}}
\singlespace
\newtheorem{thmm}{Theorem}
\newtheorem{co}[thmm]{Corollary}
\newtheorem{th}{Theorem}[section]
\newtheorem{prop}[th]{Proposition}
\newtheorem{coro}[th]{Corollary}
\newtheorem{lem}[th]{Lemma}
\newtheorem{rem}[th]{Remark}
\newtheorem{de}[th]{Definition}
\newtheorem{con}[th]{Conjecture}
\newtheorem{ex}[th]{Example}

\begin{center}
{\Large {\bf  Rational Vertex Operator Algebras and the Effective 
Central Charge}} \\
\vspace{0.5cm}

Chongying Dong\footnote{Supported by NSF grant DMS-9987656 and a research grant from the
Committee on Research, UC Santa Cruz (dong@math.ucsc.edu).}  \ \ \  and \ \ \ 
Geoffrey Mason\footnote{Supported by NSF grant DMS-9700909
and a research grant from the
Committee on Research, UC Santa Cruz (gem@cats.ucsc.edu).}\\
Mathematics Department, University of California,
Santa Cruz, CA 95064, U.S.A.

\end{center}
\hspace{1.5 cm}

\begin{abstract} 
We establish that the Lie algebra of weight one states in a 
(strongly) rational vertex operator algebra is reductive, and that 
its Lie rank $l$ is bounded above by the effective central charge 
$\c$. We show that lattice vertex operator algebras may be 
characterized by the equalities $\c= l =c$, and in particular 
holomorphic lattice theories may be characterized among all 
holomorphic vertex operator algebras by the equality $l = c$.
\end{abstract}

\section{Introduction}

The purpose of this paper is the study of certain numerical 
invariants attached to a rational conformal field theory, that is a 
(strongly) rational vertex operator algebra (RVOA). These are VOAs 
$V$ such that all (admissible) $V$-modules are completely reducible 
(see Sections 2 and 3 for a precise definition of the class of VOAs that we 
study). It is known [DLM2] that  RVOAs have only a finite number of 
simple modules $V = M^{1}, M^{2}, ... M^{r}$.  
Each has a q-graded 
character
\begin{equation}\label{1.1}
Z_{M^j}(\tau)=\tr_{M^j}q^{L(0)-c/24}=  q^{\lambda_j-c/24}\sum_{n=0}^{\infty}(\dim M^{j}_{n + 
\lambda_j})q^n
\end{equation}
for some $\lambda_j$ called the {\em conformal weight} of $M^{j}$ (where
$M^j_{\l_j}\ne 0$). The 
conformal weights are {\em a priori} arbitrary complex numbers, but 
results in [DLM3] establish that they are in fact rational numbers. 
Let $\lambda_{min}$ denote the minimum of the conformal weights and 
define
\begin{equation}\label{1.2}
\c  =  c - 24 \lambda_{min},\ \  \ll_{j}  = 
\la_j - \lambda_{min}
\end{equation}
where $c$ is the central charge of $V$. $\c$ is sometimes called the 
{\em effective central charge} in the physics literature (eg. [GN]). It is 
also proved in [DLM3] that $c$ is rational, so that the effective 
central charge of an RVOA is a rational number. We will establish a 
basic inequality satisfied by $\c$, namely that it is non-negative. 
In fact, as long as $V$ is non-trivial, that is, it is not 
finite-dimensional, we will prove that
\begin{equation}\label{1.3}
 \c > 0.
\end{equation}

The non-negativity of $\c$ is actually a consequence of more 
qualitative results which constitute our first main theorems. To 
describe them, recall [FHL] that the weight $1$ subspace $V_1$ of $V$ 
carries a natural Lie algebra structure as well as a symmetric, 
invariant bilinear form. Conversely, it is well-known (eg. [L], [Li3])
that there is a natural vertex algebra structure associated to a 
Lie algebra $\g$ and symmetric, invariant bilinear form and that this 
is frequently a VOA (ie., it also has a Virasoro vector, etc.). There 
is no restriction on the algebraic structure of $V_1 = \g$ beyond 
these conditions (we discuss an example in Section 2 in which $V_1$ 
is nilpotent of class 2). On the other hand, we will establish
\begin{thmm}\label{t1}
Let $V$ be a strongly rational vertex operator algebra. 
Then the Lie algebra $V_1$ is reductive. Moreover, any $V$-module
is a completely reducible $V_1$-module.
\end{thmm}

Granted this result, we may define the {\em Lie rank} of $V$ to be the Lie 
rank of $V_1$, i.e., the dimension of a maximal abelian subalgebra of 
$V_1$, which is of course an invariant of $V_1$. We generally denote 
the Lie rank by $l$.

\begin{thmm}\label{t2} Let $V$ be a strongly rational vertex operator algebra. Then
 $l\leq \c$.
\end{thmm}
             
Note that the non-negativity of $\c$ follows from Theorem 2.

A systematic study of the relation between the invariants $c$, $\c$ 
and $l$ may give new insights into the structure and classification 
of RCFTs and RVOAs. In particular, the method that we use to prove 
Theorem \ref{t2} can also be employed to study certain extremal situations
including the case of equality in Theorem \ref{t2}. 
Here is an example:

\begin{thmm}\label{t3} let $V$ be a strongly rational vertex operator algebra. 
Then the following are equivalent:

             (i) $\c = l = c$.
            
(ii) There is a positive-definite even lattice $L$ such 
that $V$ is isomorphic to $V_L$.
\end{thmm}

We point out that the condition $l=c$ is not sufficient
to characterize the lattice vertex operator algebras. In example
(f) in Section 4 we discuss a vertex operator algebra which satisfies 
$c=l=0$ but which is not a lattice vertex operator algebra. On the other hand,
we expect that lattice vertex operator algebras are characterized
by the condition $\c=l.$ 

Next recall that a rational $V$ is called {\em holomorphic} in case the 
adjoint module is the unique simple $V$-module. Lattice theories 
$V_L$ are holomorphic precisely when $L$ is unimodular ([D1], 
[DLM2]). For holomorphic vertex operator algebras, we always have $c 
=\c$ and $V$ is strongly rational precisely when it satisfies the 
qualitatively weaker CFT-type condition (cf. Section 2 for more 
details). Thus:

\begin{co}\label{c5} Suppose that $V$ is a holomorphic 
vertex operator algebra of CFT type. Then the following conditions are equivalent:

         (a) $l=c$.

         (b) $V$ is isomorphic to a lattice theory $V_L$ for some 
positive-definite, even, unimodular lattice $L$.
\end{co}

It is worth comparing Corollary \ref{c5} 
with the work of Schellekens [Sch], who discusses the possible 
holomorphic $c=24$ theories. In particular, the theories on Schellekens' list 
which have $l = 24$ are precisely the lattice theories $V_L$ where 
$L$ is now a Niemeier lattice (including the Leech Lattice). Thus the 
special case $c = 24$ of Corollary  \ref{c5} confirms with mathematical rigor 
the completeness of this part of Schellekens' list. See [DM2] for
further developments of this idea.

Although there is no mention of modular-invariance in any of the 
Theorems we have stated, it is the driving force behind the proofs 
each of them. While it is an important conjecture that for RVOAs the 
graded characters (\ref{1.1}) are modular functions on a congruence 
subgroup of 
$SL(2,\Z),$ this remains unknown. Our methods offer a way of 
circumnavigating this lacuna in our knowledge in some circumstances. 
A careful study of what is known about the invariance properties of 
the characters (\ref{1.1}) reveals that they constitute a {\em vector-valued} 
modular form. It transpires that this is a good substitute for 
modularity of characters in some circumstances. In particular, a 
recent result [KM] establishing polynomial growth of the Fourier 
coefficients of holomorphic vector-valued modular forms figures in 
the proofs.

The paper is organized as follows: in Section 2 we discuss background
and nomenclature for vertex operator algebras, and in Section 3 we
discuss the connections between rational vertex operator algebras,
modular-invariance, and vector-valued modular forms. Section 4 is
devoted to the proofs of Theorems \ref{t1} and \ref{t2}.  In Section 
5 we give a result on uniqueness of simple current extensions of independent 
interest and which is used in Section 6 to prove Theorem \ref{t3}.

\section{Vertex Operator Algebras of CFT type}
\setcounter{equation}{0}

First we recall from [DLMM] the following definition: a vertex 
operator algebra $V$ is said to be of {\em CFT type} in case the 
$L(0)$-grading on $V$ has no non-negative weights, and if the degree 
zero homogeneous subspace $V_0$ is $1$-dimensional. In this case, we 
have
\begin{equation}\label{2.1}
V = \bigoplus_{n=0}^{\infty}V_n
\end{equation} 
moreover $V_0$ is spanned by the vacuum vector $\1$. If $V$ is a 
vertex operator algebra of CFT type, we say that $V$ is of {\em strong} CFT 
type in case it satisfies the further condition that the $L(1)$ 
operator annihilates the homogeneous subspace $V_1$, that is
\begin{equation}\label{2.2}
L(1)V_1  =  0. 
\end{equation}

In general, one knows that $L(1)$ maps $V_1$ into $V_0$. Thus for 
vertex operator algebras of strong CFT type we have
\begin{equation}\label{2.3}
\dim  V_0/ L(1)V_1  = 1. 
\end{equation}
Li has shown [Li1] that  a simple vertex operator algebra which 
satisfies (\ref{2.3}) has, up to multiplication by a non-zero scalar, a 
unique non-degenerate invariant bilinear form $\<  ,  \>$  in the sense 
of [FLH, Section 5.3] and [Li1]. If we normalize in such a way that 
the vacuum satisfies  $\<1 , 1\> = - 1$ then in fact we have
\begin{equation}\label{2.4}
\<a , b\>=-\Res_z( z^{-1}Y(e^{zL(1)}(-z^{-2})^{L(0)}a,z^{-1})b
\end{equation}
for elements $a,b$ in $V$. It is well-known (loc. cit.) that (\ref{2.4}) 
implies that the homogeneous spaces $V_n$ and $V_m$ are orthogonal if 
$n$ and $m$ are distinct, so that the restriction of $\<  ,  \>$ to each 
$V_n$ is non-degenerate. Note that  $\<  ,  \>$  is necessarily symmetric by 
Theorem 5.3.6 of [FHL]. Furthermore, from the discussion in Section 
5.3 [loc. cit.], one sees that a vertex operator algebra $V$ of CFT 
type is of strong CFT type precisely when $V$ is self-dual in the 
sense that the contragredient module $V'$ is isomorphic to the 
adjoint module $V$.

It is easy to see that if $V$ is holomorphic  of CFT type then it is 
of strong CFT type. Indeed, the uniqueness of the simple module for 
$V$ entails that the contragredient module $V'$, which is itself 
simple, is isomorphic to the adjoint module. Now the assertion 
follows from the preceding remarks.

Let us return to the consideration of simple vertex operator algebras 
$V$ of strong CFT type. For a state $v$ in $V_n$ we define $o(v)$ to 
be the n-1th. component operator $v_{n-1}$ of $v$. Although it is 
something of a misnomer, we often refer to $o(v)$ as the {\em zero mode} of 
$v$. Note that zero modes induce a linear map on each homogeneous 
space of $V$. It is well-known that the states of weight 1 close on a 
Lie algebra under the bracket operation $[a,b] =  a_0b = o(a)b$. We 
continue to denote this Lie algebra by $V_1$. Each homogeneous space 
$V_n$ becomes a module over $V_1$ via the action of the zero mode 
i.e., $a\cdot v = o(a)v$ for $a$ in $V_1$  and $v$ in $V_n$. Thus in its 
action on $V_1$, the operator $o(v)$ for $v \in V_1$ is just $ad v$. 
Furthermore, one can check from (\ref{2.4}) that the restriction of the 
form $\<  ,  \>$ to each $V_n$ is invariant under the action of $V_1$. 
That is, we have
\begin{equation}\label{2.5}
\< o(a)u ,  v \>  =  - \< 
u, o(a)v \>
\end{equation}
for $a$ in $V_1$ and $u,v$ in $V_n$.
                    
In particular, the restriction of $\<  ,  \>$ to $V_1$ endows $V_1$ with 
a non-degenerate, symmetric, invariant, bilinear form. Indeed, it is 
well-known, and follows from (\ref{2.4}), that we have
\begin{equation}\label{2.6}                                              
\< a, b \>  = a_1b
\end{equation}
for $a, b$ in $V_1$.

\begin{ex}\label{e4.2}
{\rm We discuss an example of a vertex operator algebra of 
strong CFT type for which $V_1$ is not reductive. The reader may 
construct many other similar examples. Take a 6-dimensional complex 
Lie algebra $\g$ which is nilpotent of class 2 and such that the 
center $Z(\g)$ and the commutator subalgebra $[\g,\g]$ coincide and have 
dimension 3. Thus we may take $\g$ to be generated by elements $A$, 
$B$, $C$, with $Z(\g)$ having as basis the commutators $[A,B], [A,C], 
[B,C]$. We define a non-degenerate, symmetric, invariant bilinear form 
on $\g$ by setting
$$\<A, [B,C]\>  =  \<B, [C,A]\>  =  \<C, [A,B]\>  =  1$$
and with all other basis vectors of $\g$ being orthogonal.

Let $\hat{\frak g}={\frak g}\otimes\C[t,t^{-1}]+\C K$ be the affine
Lie algebra associated to ${\frak g}$ and 
$V(1,\C)=U(\frak g)\otimes_{{\frak g}\otimes \C[t]+\C K}\C$ the associated
Verma module where ${\frak g}\otimes \C[t]$ acts on $\C$ trivially and
$K$ acts as 1. Let $V$ be the irreducible quotient of $V(1,\C).$ Then 
$V$ is a simple vertex algebra of level 1 and $V_1$ coincides with $\g$. 
We claim that there is a Virasoro element $\omega$ in $V$ such that $V$ is a vertex 
operator algebra of strong CFT type and $c = 6$. Indeed, we may take 
$\omega$ to be the same as the usual Sugawara form for the case that 
the corresponding Lie algebra $\g$ is semi-simple, i.e.,
\begin{equation}\label{2.7}
\omega  =   \frac{1}{2} \sum_{i=1}^6u_i(-1)u^i(-1)
\end{equation}
where $u_i$ ranges over a basis of $\frak g$ and $u^i$ ranges over the dual 
basis with respect to $\<  ,  \>.$ 

One may check these assertions by an explicit but lengthy 
calculation, or by using some results in [L]. Note (loc. cit.) that 
we may modify the level to be any non-zero scalar, and also modify 
the Virasoro element so that $V$ becomes a VOA of CFT type, but not 
strong CFT type.}
\end{ex}

\section{$C_2$-rational vertex operator algebras and modular forms}
\setcounter{equation}{0}

We recall some definitions from [DLM2] and [Z]. The vertex operator 
algebra $V$ is $C_2$-{\em cofinite} in case the subspace $C_2(V)$ of $V$ 
spanned by all elements of the type $a_{-2}b$ for $a$, $b$ in $V$ has 
finite codimension in $V$. $V$ is {\em rational} if any admissible 
module is a direct sum of irreducible admissible modules. 
Let us introduce the term {\em strongly 
rational} for a vertex operator algebra $V$ which satisfies the 
following conditions:

1. $V$ is of strong CFT type.
                       
2. $V$ is $C_2$-cofinite. 

3. $V$ is rational.

It is likely that these conditions are not independent. For example, 
it may well be that conditions 2 and 3 are equivalent. We will see 
that the class of vertex operator algebras which are strongly 
rational in the above sense have a number of good properties which 
can be exploited.

Let us fix for the duration of this section a strongly rational VOA 
$V$ with central charge $c$.  Let notation for simple $V$-modules be 
as in (\ref{1.1}). The basic results concerning the modular-invariance 
properties of the characters (\ref{1.1}) and related trace functions were 
established in [Z], with some refinements incorporated in [DLM3]. To 
describe these, let $u$ be a state in $V$ with vertex operator
\begin{equation}\label{3.1}
Y(u,z) = \sum_{n \in \Z} u_nz^{-n-1}.
\end{equation}

Define the zero mode $o(u)$ to be the sum of the zero modes of the 
homogeneous components of $u$, which were defined in Section 2. We 
can then define the 1-point function $Z_{M^j}(u,\tau)$ associated to $M_j$ 
as follows:
\begin{equation}\label{3.2}
Z_{M^j}(u, \tau)  =  \tr_{M^j}o(u) q^{L(0) - c/24}=
q^{ - c/24 + \lambda_j} \sum_{n\geq 0}\tr_{M^j_{n+ \lambda_j}}o(u)q^n.
\end{equation}
In case the state $u$ is the vacuum vector $1$, we have
\begin{equation}\label{3.3}
Z_{M^j}(\tau)  =Z_{M^j}(\1, \tau).
\end{equation}
Here and below, $\tau$ denotes an element in the complex upper 
half-plane $\H$ and $q = e^{2 \pi i \tau}$. 

In order to discuss the modular properties of $Z_{M}(v,\tau)$ we briefly recall the genus
one vertex operator algebra $(V,Y[\ ],\1, \omega-c/24)$ from [Z]. 
The new vertex 
operator associated to  a homogeneous element $a$ is given by
\[ Y[a,z] = \sum_{n\in\Z}a[n]z^{-n-1} = Y(a,e^{z} -1)e^{z\wt{a}}
\] while the Virasoro element is $\tilde{\omega} = \omega-c/24$.
Thus 
\[
a[m] = \Res_z\left(Y(a,z)(\ln{(1+z)})^m(1+z)^{\wt{a}-1}\right)
\]
and \[
a[m] = \sum_{i=m}^\infty c(\wt{a},i,m)a(i)
\]
for some scalars $c(\wt{a},i,m)$ such that $c(\wt{a},m,m)=1.$ 
In particular,
\[
a[0]=\sum_{i\geq 0}{\wt{a}-1 \choose i}a(i).
\]
We also write
\[
L[z] = Y[\omega,z] = \sum_{n\in\Z} L[n]z^{-n-2}.
\]
Then the $L[n]$ again generate a copy of the Virasoro algebra with
the same central charge $c.$ Now $V$ is graded by
the $L[0]$-eigenvalues, that is
\[
V=\bigoplus_{n\in\Z}V_{[n]}
\]
where $V_{[n]}=\{v\in V|L[0]v=nv\}.$ We also write ${\wt}[a]=n$
if $a\in V_{[n]}.$ It should be pointed out that for any
$n\in\Z$ we have 

\[
\sum_{m\leq n}V_n=\sum_{m\leq n}V_{[n]}.
\] 

It is established in 
[DLM3] and [Z] that the functions (\ref{3.2}) are holomorphic in $\H$ 
and that the following holds: for any state $u$, homogeneous of 
weight $k$ with respect to the square bracket Virasoro operator 
$L[0]$,  
and for any matrix $\gamma  =\left(\begin{array}{cc}
 a & b\\
c & d
\end{array}
\right)$ in the modular group $SL(2,\Z)$, there are scalars
$\rho_{i,j}(\gamma)$, $1\leq i,j \leq r$ independent of $u$ and
$\tau,$ and an equality
\begin{equation}\label{3.4}
Z_{M^i}(u, \frac{a\tau + b}{c\tau+d})=(c\tau + d)^k\sum_{j=1}^r
\rho_{i,j}(\gamma)Z_{M^j}( u, \tau).
\end{equation}

It will be worthwhile to state these results in 
another form; we refer the reader to [KM] for relevant background. 
Let $Hol$ and $Hol_q$ denote respectively the space of holomorphic 
functions on $\H$ and the subspace of functions which are in 
addition 'meromorphic at infinity' i.e., have a q-expansion
\begin{equation}\label{3.5}
\sum_{n\geq b} a_nq^{n/N}
\end{equation}
for some $b \in \Z$ and $N \in \N$. It is well-known that $Hol$ is a 
right $\Gamma$-module with respect to the kth. slash operator
$$f|_k \gamma (\tau) =  (c\tau + 
d)^{-k} f(\gamma\tau).$$

Let $Hol_k$ denote the space $Hol$ considered as a right 
$\Gamma$-module in this way. What (\ref{3.4}) says is that for each state $u 
\in V$ satisfying $L[0]u = ku$, the 1-point functions $Z_{M^i}(u, \tau)$ 
span a finite-dimensional $\Gamma$-submodule $R_u$ of $Hol_k$. 
Moreover $R_u \subset Hol_q$. (Note that $Hol_q$ is not a $\Gamma$-submodule 
of $Hol_k$.) In fact, the matrices $(\rho_{i,j}(\gamma))$ may be chosen 
so that they define a representation of $\Gamma$ on the so-called 
conformal block $B$ (cf. [Z]), and a choice of $u$ as above 
determines a homomorphism of $\Gamma$-modules from $B$ to $R_u$.

Yet another way to state the above results is in terms of so-called 
(meromorphic) {\em vector-valued modular forms} of weight $k$: by this we 
mean a tuple of functions
$${\F}(\tau) =  (F_1(\tau), ... , F_r(\tau))$$
 together with a 
representation $\rho: \Gamma \to GL(r,\C)$ such that each component 
function $F_i(\tau)$ lies in $Hol_q$ and the following compatibility 
condition holds:
\begin{equation}\label{3.6}
\F^T |_k \gamma(\tau) =\rho(\gamma)\F^T(\tau)
\end{equation}
for all $\gamma \in \Gamma$. $T$ denotes transpose of vectors, and the 
slash operator acts on vectors of functions componentwise. We call 
the vector-valued modular form {\em holomorphic}, or {\em cuspidal}, if the 
component functions $F_i(\tau)$ are holomorphic at infinity or vanish 
at infinity respectively i.e., the integer $b$ in (\ref{3.5}) is 
non-negative or positive respectively. We note that we may take the 
weight $k$ here to be any rational number.
We can now restate (\ref{3.4}) in the following way:

\begin{prop}\label{p3.1}
For each state $u \in V$ which is homogeneous of 
weight $k$ with respect to the operator $L[0]$, the $r$-tuple 
${\bf Z}(u,\tau) = (Z_{M^1}(u,\tau), ..., Z_{M^r}(u, \tau))$ is a vector-valued 
modular form of weight $k$ with respect to the representation $\rho$.
\end{prop}

It is known (cf. [KM]) that all modular forms in the usual sense 
give rise to vector-valued modular forms on a subgroup of $\Gamma$ of 
finite index, however there are vector-valued 
modular forms which do not arise in this way. Vector-valued modular 
forms nevertheless enjoy some of the properties of ordinary modular 
forms. The following result, which extends a classical result of 
Hecke for ordinary modular forms and which plays a role in the 
present paper, illustrates this idea.

\begin{prop}\label{p3.2} Suppose that $F(\tau)$ is a component of a 
holomorphic vector-valued modular form of weight $k.$ Then the Fourier 
coefficients $a_n$ of $F(\tau)$ satisfy the growth condition $a_n = O(n^{\alpha})$ 
for a constant $\alpha$ independent of $n$.
\end{prop}

\pf  This is the main result of [KM]. \qed

Compare this to the growth of the Fourier coefficients of the 
meromorphic modular form $\eta(\tau)^{-1}$. Here, $\eta(\tau)$ is the 
Dedekind eta function
$$\eta(\tau) = q^{1/24}\phi(q) = q^{1/24}\prod_{n=1}^{\infty}(1 - q^n)$$
 and
$$\phi(q)^{-1}  =\sum_{n=0}^{\infty}p(n) q^n$$
with $p(n)$ the usual unrestricted partition function. It is 
well-known(cf. [K]) that the growth of $p(n)$ is exponential in $n$.

Propositions \ref{p3.1} and \ref{p3.2} do not necessarily apply to the 
vector-valued modular forms ${\bf Z}(u,\tau)$ {\em per se:} they generally have 
poles. However, we can multiply by a suitable cusp-form (say, a power 
of the eta function) in order to remove the poles. What we will 
actually be doing is very close to this, but in the guise of a 
Lie-theoretic argument involving Heisenberg and Virasoro Lie algebras.

Because of this circumstance, it will be useful to look more closely 
at the poles of the $q$-expansions of the functions $Z_{M^i}(\tau)$. 
Indeed, from (\ref{3.3}) and (\ref{1.2}) we have
$$Z_{M^i}(\tau)  =  q^{ - \c/ 24 + \ll_i }\sum_{n\geq 0}(\dim M^i_{n + \lambda_i})q^{n}.$$

Since $\ll_i \geq 0,$ it follows that the pole at infinity for 
$Z_{M^i}(\tau)$ has order no worse than $\c/ 24.$ 
Alternatively, we can 
assert that $\eta(\tau)^{\c} Z_{M^i}(\tau)$ is holomorphic in ${\bf H}\cup 
\{i\infty\}.$ As a result of the transformation law for $\eta(\tau)$, it 
follows that $\eta(\tau)^{\c}{\bf Z}(u,\tau)$ is a holomorphic vector-valued 
modular form of weight $k + \c/2$. By Proposition \ref{p3.2} we can conclude

\begin{coro}\label{c3.3} The Fourier coefficients of 
$\eta(\tau)^{\c}Z_{M^i}(u,\tau)$ satisfy a polynomial growth condition 
$a_n = O(n^{\alpha})$.
\end{coro}

There is a related result ([KM, Lemma 4.2]) that we use:

\begin{lem}\label{l3.4} Let $\F$ be a vector-valued modular form such that the 
Fourier coefficients of each component function $F_i(\tau)$ are 
non-negative real numbers. If $\F$ is holomorphic of weight $0$ 
then it is a constant, and if $\F$ is cuspidal then $\F = 0$.
\end{lem}

In a somewhat different direction, we recall an identity due to Zhu 
[Z]. It may be construed as a VOA-theoretic analog of the Killing 
form familiar in Lie theory. Namely, for states $u, v$ in $V$ we have
\begin{equation}\label{3.7}
\tr_{M^i} o(u)o(v)q^{L(0) - c/24}=Z_{M^i}(u[-1]v, \tau) - \sum_{k\geq 1}E_{2k}(\tau)Z_{M^i}(u[2k-1]v, \tau).
\end{equation}

In (\ref{3.7}), the meaning of the operators $u[k]v$ has already been 
described. The functions $E_{2k}(\tau)$ are the usual Eisenstein series 
of weight $2k,$ normalized as in [DLM3]:
\[
E_{2k}(\tau) =\frac{-B_{2k}}{2k!}+\frac{2}{(2k-1)!}
\sum_{n=1}^\infty \sigma_{2k-1}(n)q^n
\]
where $\sigma_k(n)$ is the sum of the $k$-powers of the divisors of
$n$ and $B_{2k}$ a Bernoulli number.
Of course, $E_{2k}(\tau)$ is a 
holomorphic modular form on $SL(2,\Z)$ if $k>2$. However - and this 
will be important for us - $E_2(\tau)$ is not modular. Its 
transformation with respect to the $S$ matrix is as follows:
\begin{equation}\label{3.8}
E_2(-1/\tau)  =  \tau^2E_2(\tau)- \frac{\tau}{2 \pi i}.
\end{equation}

\section{Strongly rational vertex operator algebras and reductive 
Lie algebras}
\setcounter{equation}{0}

The main purpose of this section is to prove Theorems \ref{t1} and \ref{t2}. 
Concerning Theorem \ref{t1}, we have to show that the nilpotent radical 
${\cal N}$ of the Lie algebra $V_1$ is trivial. We use basic facts about 
$V_1$, as discussed in Section 2, without further comment. Proceeding 
by contradiction, pick a non-zero element $u$ in ${\cal N}$. Then for any 
other element $v$ in $V_1$ we see that the trace of the operator 
$o(u)o(v)$ acting on each of the homogeneous spaces $V_n$ is $0$. In 
other words, the left-hand-side of equation (\ref{3.7}) is zero. Moreover, if 
$k>1$ is an integer then because $V$ is of strong CFT type, the 
element $u[2k-1]v$ has $L[0]$-weight $2-2k < 0$ and hence is $0$. In 
addition, the element $u[1]v$ coincides with $u_1v  =  \<u , v\>$, the 
latter equality by (\ref{2.6}). As $\<  ,  \>$ is non-degenerate, we can 
choose $v$ so that  $\< u , v \> = 1$. With this choice of $v$, eqn. 
(\ref{3.7}) simplifies to read
\begin{equation}\label{4.1}
Z_{M^i}(u[-1]v, \tau)  =  E_2(\tau) Z_{M^i}(\tau).
\end{equation}

The point is that because $E_2(\tau)$ enjoys an exceptional
transformation law (\ref{3.8}), it cannot participate in an equation
of the type (\ref{4.1}). This is how we get the desired
contradiction. In detail: take any index $i$, replace the state $u$ in
(\ref{3.4}) by $u[-1]v$ with $u$, $v$ as above, so that $u[-1]v$ has
$L[0]$ weight 2, and take $\gamma$ to be the $S$ matrix. Several
applications of (\ref{3.4}) and (\ref{4.1}) yield
\begin{eqnarray*}
& &\tau^2 \sum_{j=1}^r S_{i,j}Z_{M^j}(u[-1]v, \tau)=Z_{M^i}(u[-1]v, -1/\tau)\\
& &\ \ \ \ \ \  =E_2(-1/\tau) Z_{M^i}(-1/\tau)\\
& &\ \ \ \ \ \ =(\tau^2 E_2(\tau)  -  \frac{\tau}{2 \pi i})\sum_{j=1}^rS_{i,j}Z_{M^j}(\tau)\\
& &\ \ \ \ \ \ =\tau^2\sum_{j=1}^r S_{i,j}Z_{M^j}(u[-1]v, \tau)
                                - \frac{\tau}{2 \pi i}\sum_{j=1}^r S_{i,j}Z_{M^j}(\tau).
\end{eqnarray*}
It follows that in fact
\begin{equation}\label{ae}
\sum_{j=1}^r S_{i,j}Z_{M^j}(\tau)  =0.
\end{equation}
However the left-hand-side of (\ref{ae}) is equal to 
$Z_{M^i}(-1/\tau)$, and this is clearly not zero. 
Thus $V_1$ is reductive. 

We next prove that any $V$-module $M$ is a completely reducible $V_1$-module. 
Since $V$ is rational we can assume that $M$ is irreducible.
Let $\h$ be a maximal abelian
subalgebra, so that $\h$ has dimension $l$. It is enough to prove that
$M$ is a completely reducible $\h$-module. It is well-known that
the restriction of the non-degenerate form $\< , \>$ to $\h$ is also
non-degenerate. As a result, the component operators of the vertex
operators $Y(u,z)$ on $M$ for $u$ in $\h$ close on an affine Lie algebra
$\hat\h$. That is, they satisfy the relations
$$[u_m, v_n]  =  m\delta_{m,-n}\<u , v\>.$$

Essentially by the Stone-von-Neumann theorem, there is a tensor 
decomposition of $V$,
\begin{equation}\label{4.2}
M =  M(1)\otimes_{\C} \Omega_M
\end{equation} 
where $M(1)=\C[u_{-m}|u\in {\frak h},m>0]$ is the Heisenberg vertex
operator algebra of rank $l$ generated by $\h$ and $\Omega_M$ is the
so-called {\em vacuum space} consisting of those states $v \in M$ such that
$u_nv = 0$ for all $u \in \h$ and all $n>0$. See [FLM, theorem 1.7.3]
for more details. For $\alpha\in\h$ let $M(\alpha)$ be the
generalized eigenspace for $\h$ with eigenvalue $\alpha.$ 
That is, $v\in M(\alpha)$ if and only if
there exists a positive integer $m$ such that $(u_0-\<h,\alpha\>)^mv=0$
for all $u\in\h.$ Then $M(\alpha)=M(1)\otimes \Omega_M(\alpha)$ where
$\Omega_M(\alpha)=M\cap \Omega_M.$ Let $L_M$ be the subset of $\h$ consisting
of $\alpha$ such that $M(\alpha)\ne 0.$ We set $L=L_V.$ Then
$$M=\bigoplus_{\alpha\in L_M}M(1)\otimes \Omega_M(\alpha)$$
and $u_mM(\beta)\subset M(\alpha+\beta)$ for $\alpha\in L,$ $\beta\in L_M,$
$u\in V(\alpha)$ and $m\in \Z.$  

First we take $M=V.$ Note that $M(1)$ is a vertex operator algebra
with Virasoro vector $\omega'=\frac{1}{2}\sum_{i=1}^l(u^i_{-1})^2\1$ where $\{u^1,...,u^l\}$ is
an orthonormal basis of $\h$ with respect to $\<,\>.$ 
Set $\omega''=\omega-\omega'.$ We will use the  notation
$$Y(\omega', z)  = \sum_{n\in\Z}L'(n) z^{-n-2},
Y(\omega'', z)  = \sum_{n\in\Z}L''(n) z^{-n-2}.$$
Then $\omega''$ is a Virasoro vector with central charge $c-l$
and $\Omega_V(0)$ is a vertex operator algebra with Virasoro vector
$\omega''.$ 

Fix $\beta\in L.$ Since $V$ is simple,  for any nonzero $w\in M(\beta),$
$V$ is spanned 
by $u_qw$ for $u\in V$  and $q\in \Z$ by 
Corollary 4.2 of [DM1] or Proposition 4.1 of [L2]. Thus
$V(-\beta)$ must be nonzero, otherwise $V(0)$ would be zero.
By Proposition 11.9 of [DL], for any $\alpha\in L$ and
$0\ne u\in V(\alpha)$ there exists $m\in \Z$
such that $0\ne u_mw \in V(\alpha+\beta).$ Thus $V(0)=M(1)\otimes \Omega_V(0)$
is a simple vertex operator subalgebra of $V$ and each $V(\alpha)
=M(1)\otimes \Omega_V(\alpha)$ is an irreducible
$V(0)$-module.  More precisely, as a $V(0)$-module,
\begin{equation}\label{-1}
V(\alpha)=M(1,\alpha)\otimes \Omega_V(\alpha)
\end{equation}
where $M(1,\alpha)=\C[u_{-n}|{u\in \h}, n>0]e^{\alpha}$ is
the highest weight 
irreducible $M(1)$-module such that $h\in \h$ acts on $e^{\alpha}$ 
as $\<h,\alpha\>$ and $\Omega_V(\alpha)$ is an irreducible 
$\Omega_V(0)$-module. Thus $V$ is completely reducible $\h$-module.  

Applying the same argument to any irreducible $V$-module $M$ shows
that $M$ is a completely reducible $\h$-module. This
finishes the proof of Theorem \ref{t1}. \qed

\begin{rem}\label{r4.1} {\rm From the proof of Theorem \ref{t1}, we see that 
$L$ is an additive subgroup of $\h$ and each $L_M$ is a coset of $L$ in
$\h.$}
\end{rem}  

We turn next to the proof of Theorem \ref{t2}. Note that both factors on the right-hand-side of
(\ref{4.2}) in the case $M=V$ 
are invariant under the $L(0)$ operator, so that (\ref{4.2}) is an
equality of $q$-graded spaces. Because the $q$-graded character of $M(1)$
is equal to $\phi(q)^{-l}$, it follows from (\ref{4.2}) that
\begin{equation}\label{4.3}
\eta(\tau)^{\c}Z_V(\tau) =q^{(l-c)/24}\eta(\tau)^{\c -l} ch_q\Omega_V
\end{equation}
where we have used $ch_q(\Omega_V)$ for the $q$-graded character of 
$\Omega_V$. Now according to Corollary \ref{c3.3} the Fourier coefficients of 
the l.h.s. of (\ref{4.3}) have polynomial growth, so the same must be true 
of the Fourier coefficients on the r.h.s. But $\eta(\tau)^{-1}$ has 
exponential growth of Fourier coefficients, as indeed does 
$\eta(\tau)^s$ whenever
$s < 0$. We conclude that $\c - l \geq 0$, which completes the proof of 
Theorem \ref{t2}.

A second way to complete the proof of the theorem runs as 
follows: (\ref{4.2}) and (\ref{4.3}) apply not only to $V$ but to each simple 
$V$-module. As a result, the (non-zero) vector-valued modular form 
$\eta(\tau)^{\c}{\bf Z}(\1,\tau)$ has components whose Fourier coefficients 
are non-negative integers and which, assuming that $\c < l$, is 
cuspidal. This contradicts Lemma \ref{l3.4}. Furthermore, if $\c = 0$ the same argument shows that 
${\bf Z}(\1,\tau)$ is necessarily constant, so that $V$ is $1$-dimensional. 
As a result, (\ref{1.3}) holds if $V$ is non-trivial.  \qed

\bigskip

We discuss the invariants $c, \c, l$ for some standard VOAs.

(a) If $V$ is strongly rational and holomorphic, Theorem \ref{t2} tell us 
that  $l\leq c = \c$.

\bigskip

(b) Suppose that $V$ is a unitary VOA (cf. [Ho] for background on 
such theories). Then one knows (loc. cit.) that $c$ and all conformal 
weights $\lambda_i$ are non-negative. In particular, as in the 
holomorphic case, we have $l\leq c = \c$ for unitary theories.

\bigskip

(c) Suppose that $V = V_L$  is the vertex operator algebra associated 
to a positive-definite, even lattice $L$. Then $V$ is strongly 
rational and $l = c = \c= rank L$. See [D1], [DLM1].

\bigskip

(d) Consider the simple vertex operator algebra  $V$ = $L({\frak g}, k)$ 
associated to a finite-dimensional simple Lie algebra $\frak g$ with level 
$k$. Then $l$ is the Lie rank of $\g$ and $c = k \dim\g / (k + h^{\vee})$ 
as long as $k + h^{\vee}$  is non-zero. Here, $h^{\vee}$  is the dual Coxeter 
number of $\g$. $V$  is rational (and hence strongly rational) 
precisely when $k$ is a nonnegative integer (cf. [FZ], [DLM1]).

\bigskip

(e) Consider the Virasoro vertex operator 
algebra $V=L(c,0)$ with $c= c_{p,q} = 1 - 6(p - q)^2/pq$ for coprime
positive integers $p, q$ with $1<p<q$, say. (See Section 6 for the definition
of $L(c,0).$) It is known ([W], [DLM1])
that $V$ is strongly rational. The conformal weights of the simple
$V$-modules are the rational numbers $\lambda_{m,n} = ((np - mq)^2 -
(p - q)^2)/ 4pq$ for integers $m, n$ in the range $1\leq m\leq p - 1,
1\leq n\leq q- 1$. These theories are unitary precisely when $q - p =
1$, in which case $\c = 1 - 6/ pq$. We claim that this is the
effective central charge in all cases. This is equivalent to the
equality $\lambda_{min} = (1 - (p - q)^2)/ 4pq$, which in turn is
equivalent to the assertion that we may choose integers $m, n$ in the
indicated range such that $|np - mq| = 1$. This latter fact is
well-known.

\bigskip

(f) Consider the vertex operator algebra $V=L(c_{2,5},0)^{\otimes 60}\otimes
(\V)^{\otimes 11}$  where $\V$ is the moonshine 
vertex operator algebra constructed in [FLM]. Since the central charges
of $L(c_{2,5},0)$ and $\V$ are $-22/5$ and $24$ respectively, $V$ has
zero central charge. Note that both  $L(c_{2,5},0)$ and $\V$
are strongly rational with weight one subspaces being zero (cf.
[W] and [D2]. Thus $V$ is a strongly rational vertex operator algebra
with zero weight one subspace. As a result the Lie rank $l$ of $V$ is
zero. Since the Lie rank of any lattice vertex operator
algebra $V_L$ is equal to the rank of $L$ we conclude that $V$ 
is not a lattice vertex operator algebra. The effective central charge
$\c$ is equal to $528.$ One can get a strongly rational vertex
operator algebra with $c=l=0$ and $\c$ equal to any positive 
multiple of $528$ by tensoring $V$ with itself.

\bigskip 

We learn from these examples that one cannot expect to prove a 
universal inequality  between $l$ and $c$. Of course we always have 
$c\leq \c$.

\section{Uniqueness of simple current extensions}
\setcounter{equation}{0}

We 
establish a general result of independent interests which will be used 
in the proof of Theorem \ref{t3}. 
In the language of physics, the result to be proved in this section 
establishes the uniqueness of vertex operator algebras defined by 
simple current extensions.

We first recall from [FHL] the definition of intertwining operators.
 \begin{de} 
Let $V$ be a vertex operator algebra and let $W_{1},W_{2}$ and $W_{3}$ 
be weak $V$-modules.
An {\em intertwining operator } of type ${W_{3}\choose W_{1}W_{1}}$ 
is a linear map $I$ from $W_{1}$ to $(\Hom (W_{2},W_{3}))\{z\}$
satisfying the following properties:
for $w_{1}\in W_{1},\; w_{2}\in W_{2}$, 
\begin{eqnarray}
I(w_{1},z)w_{2}\in W_{3}\{z\}
\end{eqnarray}
for $v\in V,\; w_{1}\in W_{1}$ where $W\{z\}$ is the space
of formal power series in $z$ with complex powers
and coefficients in $W$ for a vector space $W,$
\begin{eqnarray}
[L(-1),I(w_{1},z)]={d\over dz}I(w_{1},z),\label{ederivative}
\end{eqnarray}
\begin{eqnarray}\label{ejacobiinter}
& &z^{-1}_0\delta\left(\frac{z_1-z_2}{z_0}\right)
Y(v,z_1)I(w_{1},z_2)-z^{-1}_0\delta\left(\frac{z_2-z_1}{-z_0}\right)
I(w_{1},z_2)Y(v,z_1)\nonumber\\
&=&z_2^{-1}\delta\left(\frac{z_1-z_0}{z_2}\right)
I(Y(v,z_0)w_{1},z_2).
\end{eqnarray}

\end{de}

All intertwining operators of type 
${W_{3}\choose W_{1}W_{2}}$ form a vector space denoted by
$I_{W_{1}W_{2}}^{W_{3}}$. The dimension of $I_{W_{1}W_{2}}^{W_{3}}$
is called a {\em fusion rule}, 
denoted by $N_{W_{1},W_{2}}^{W_{3}}$.
Clearly, fusion rules only depend on the equivalence class 
of each $W_{i}$.

\begin{de} A vertex operator algebra $V$ is graded by an abelian group
$G$ if $V=\oplus_{g\in G}V^g,$ and  $u_nv\in V^{g+h}$ for
any $u\in V^g, v\in V^h$ and $n\in\Z.$ 
\end{de}

It is easy to see from the definition that $V^0$ is a vertex operator
subalgebra of $V$ and each $V^g$ is a strong $V^0$-module for
$g\in G.$ Moreover, the restriction $Y_{g,h}(u,z)$ of $Y(u,z)$ to $V^h$ for
$u\in V^g$ is an intertwining operator of type ${V^{g+h}\choose V^g V^h}$ 
for $V^0$-modules $V^g,V^h,V^{g+h}.$ 

\begin{prop}\label{in} Let $V=\sum_{g\in G}V^g$ be a simple $G$-graded vertex 
operator algebra such that $V^g\ne 0$ for all $g\in G$ and 
$N_{V^g,V^h}^{V^k}=\delta_{g+h,k}$ for all $g,h,k\in G.$ Then the vertex operator
algebra structure of $V$ is determined uniquely by the $V^0$-module
structure of $V.$ That is, if $(\bar V,\bar Y,{\bf 1},\o)$ is also a simple 
vertex operator algebra with $\bar V=V$ as vectors spaces and 
$Y(u,z)=\bar Y(u,z)$ for all $u\in V^0$ then  
$(V,Y,{\bf 1},\o)$ and $(V,\bar Y,{\bf 1},\o)$ are isomorphic.
\end{prop}

\pf Note that  both $\bar Y_{g,h}$ and $Y_{g,h}$ are intertwining operators
of type  ${V^{g+h}\choose V^g V^h}$ for $V^0$-modules $V^g,V^h,V^{g+h}$
for any $g,h\in G.$ By assumption there exists constants $c_{g,h}$ such that
$\bar Y_{g,h}=c_{g,h}Y_{g,h}.$ We shall prove that $c: G\times G\to \C^{*}$
defines a normalized two cocycle on $G.$

It is clear that  $c_{0,h}=1$ for all $h\in G$ as $\bar Y(v,z)=Y(v,z)$ for
all $v\in V^0.$ By skew symmetry (see [FHL]), for $u\in V^h$ and $v\in V^g,$
$$\bar Y(u,z)v=e^{zL(-1)}\bar Y(v,-z)u=c_{g,h}e^{zL(-1)}Y(v,-z)u=c_{g,h}Y(u,z)v.$$
Thus $c_{g,h}=c_{h,g}$ for all $g,h\in G.$ In particular, 
this shows that $c_{h,0}=1.$

Using commutativity (see [FLM], [FHL] and [DL]) for $u\in V^g$ and $v\in V^H$
gives a non-negative integer $n$ such that
$$(z_1-z_2)^n\bar Y(u,z_1)\bar Y(v,z_2)=(z_1-z_2)^n\bar Y(v,z_2)\bar
Y(u,z_1)$$ 
and that
$$(z_1-z_2)^nY(u,z_1)Y(v,z_2)=(z_1-z_2)^nY(v,z_2)Y(u,z_1).$$
Applying both identities to $V^k$ yields
$$c_{g,h+k}c_{h,k}=c_{h,g+k}c_{g,k}=c_{g+k,h}c_{g,k}.$$
That is, $c$ is a 2-cocycle.

Since $\C^*$ is injective in the category of abelian groups, there
exists a function $f: G\to \C^*$ such that $f(0)=1$ and
$c_{g,h}=f(g)f(h)/f(g+h)$ for $g,h\in G.$ Define a linear map $\sigma$
from $\bar V$ to $V$ by sending $v$ to $f(g)v$ for $v\in V^g$ and
$g\in G.$ It is enough to prove that $\sigma$ is an isomorphism of
vertex operator algebras.

Clearly, $\sigma$ is a linear isomorphism. Let $u\in V^g$ and $v\in V^h.$ Then
\begin{eqnarray*}
& &\ \ \ \ \  \sigma(\bar Y(u,z)v)=f(g+h)\bar Y(u,z)v\\
& &=f(g+h)c_{g,h}Y(u,z)v\\
& &=f(g)f(h)Y(u,z)v\\
& &=Y(\sigma u,z)\sigma v.
\end{eqnarray*}
Thus $\sigma$ is an isomorphism of vertex operator algebras.
\qed

In order to discuss a consequence of Proposition \ref{in} we recall
from [B] and 
[FLM] the vertex operator algebra $V_L$ associated to an even positive 
definite lattice $L.$ So $L$ is a free abelian group of finite rank
with a positive definite $\Z$-bilinear form $(,)$ such that
$(\alpha,\alpha)\in 2\Z$ for $\alpha\in L.$ Set ${\frak h}=\C\otimes_{\Z}L$
and 
let  $\hat{\frak h}={\frak h}\otimes \C[t,t^{-1}]\oplus\C c$ the corresponding
affine Lie algebra. Let $M(1)=\C[h(-n)|h\in{\frak h}, n>0]$ be the unique
irreducible module for $\hat{\frak h}$ such that $c$ acts as $1$ and ${\frak h}\otimes t^0$ acts trivially. Then 
as a vector space,
$$V_{L}=M(1)\otimes \C[L]$$
where  $\C[L]$ is the group algebra $\C[L].$
Then $V_L$ is a strongly rational vertex operator algebra (see [B], [FLM],
[D1] and [DLM1]). 

It is a fact that $M(1)$ is a vertex operator subalgebra of $V_L$
(see [FLM]) and $V_L$ is $L$-graded such that $V_L^{\alpha}=M(1)\otimes 
e^{\alpha}$ for $\alpha\in L$ where $e^{\alpha}$ denotes the basis element
in $\C[L]$ corresponding to $\alpha.$ In fact, $V_L^{\alpha}=M(1)\otimes 
e^{\alpha}$ is an irreducible $M(1)$-module.

\begin{coro}\label{ca} Let $V$ be simple vertex operator algebra which
has a subalgebra isomorphic to $M(1),$ such that $V$ isomorphic to
$V_L$ as $M(1)$-modules for some even positive definite lattice $L.$ 
Then $V$ and $V_L$ are isomorphic vertex operator algebras.
\end{coro}

\pf By assumption, $V\equiv \sum_{\alpha\in L}V_L^{\alpha}.$ By Proposition
8.15 of [G], $N_{V_L^{\alpha},V_L^{\beta}}^{V_L^{\gamma}}=\delta_{\alpha+\beta,\gamma}.$  Proposition \ref{in} then gives the result.
\qed

\section{Proof of Theorem 3}
\setcounter{equation}{0}

In order to prove Theorem \ref{t3} we first recall some results about
highest weight modules for the Virasoro algebra. 

For a pair of complex numbers $c, h$ we denote by $Ver(c, h)$ the 
Verma module of highest weight $(c, h)$ over the Virasoro algebra 
$Vir_c$ of central charge $c$.
Thus $Ver(c,h)$ is generated as 
$Vir_c$-module by a state $v_0$ satisfying $L(n)v_0 = 0$ for $ n>0$ 
and $L(0)v_0 = hv_0$. The structure of these modules was elucidated 
in the work of Feigin-Fuchs [FF],  Rocha-Caridi and Wallach [RW1],[RW2],
and Kac and Raine [KR], and their results 
play an important role in our proof of Theorems 3 and 4. The case $h = 0$ is 
particularly important. One knows (loc. cit.) that  $Ver(c, 0)$ 
contains a singular vector $L(-1)v_0 =v_1$ at level $1$, and by 
Frenkel-Zhu [FZ] the quotient module $V(c, 0) / U(Vir_c)v_1$ 
carries the structure of a vertex operator algebra. We denote this 
vertex operator algebra by $V(c, 0)$. It contains a unique maximal 
$Vir_c$-submodule $J(c, 0)$, and the quotient module $L(c, 0) = V(c, 
0)/ J(c, 0)$ also carries the structure of vertex operator algebra. 
The rationality of $L(c, 0)$ was investigated by Wang [W]. We collect 
some of the results of these authors in the following

\begin{prop}\label{vir} One of the following holds:

(i) $J(c,0) = 0$, so that $V(c, 0) = L(c, 0)$ is a simple module over 
$Vir_c$. In this case, the $q$-graded character of $L(c,0)$ is equal 
to $(1 - q)/\phi(q)$ and the coefficients have exponential growth. 
Moreover, $L(c, 0)$ is not a rational vertex operator algebra.

(ii) $J(c, 0)\ne 0$. In this case $L(c, 0)$ is a rational vertex 
operator algebra. The central charge has the form 
$c = c_{p, q} = 1 -6(p - q)^2/pq$ for a pair of coprime integers $p, q$ 
satisfying $2\leq p<q$.
\end{prop}

In order to prove Theorem \ref{t3} we first 
assume only that $\c=l.$ We continue with the notation of Section  4.
Let $W$  denote the Virasoro vertex 
operator algebra generated by $\omega''$.  Then $W$ is a highest weight
module for the Virasoro algebra with highest weight $0.$ 
Of course it is possible 
that $\omega'' = 0$, in which case $W$ is nothing but the complex 
numbers. 

\begin{lem}\label{llast} Let $V$ be a strongly rational vertex operator
algebra satisfying $\c=l.$  Then there exist coprime positive integers
$p,q$ with $2\leq p<q $ such that $W$ is isomorphic to $L(c_{p,q},0).$
\end{lem}

\pf Recall from (\ref{4.2}) that $V=M(1)\otimes \Omega_V$ and $\Omega_V$ is a 
module for the Virasoro algebra $Vir''.$  
From the proof of Corollary \ref{c3.3} 
we know that  the Fourier coefficients
of $$\eta(\tau)^{\c}Z_V(\tau) =q^{(l-c)/24}\eta(\tau)^{\c -l} ch_q\Omega_V
=q^{(l-c)/24}\ch_q\Omega_V $$
have polynomial growth. If $c''=c-l$ is not equal to any $c_{p,q}$ for all
coprime positive integers $p,q\geq 2$ then 
$$\ch_qW=\sum_{n\geq 0}(\dim W_n)q^n=\frac{1-q}{\phi(q)}$$
by Proposition \ref{vir}. We then have that
the Fourier coefficients of $\eta(\tau)^{\c}Z_V(\tau)$ has exponential
growth. This shows that $c''=c_{p,q}$ for some coprime positive integers
$q>p>1.$ The same argument also shows that $W$ is an irreducible
highest weight module for $Vir''.$ Thus $W$ is isomorphic
to $L(c_{p,q},0).$ \qed

\bigskip

We now turn to the proof of Theorem \ref{t3}. In this case, $c''=c_{p,q}=0.$
Thus $\omega''=0$ and $W=\C.$ As a result we see that 
$$V=\sum_{\alpha\in L}M(1,\alpha)=\oplus_{\alpha\in L}M(1)\otimes \Omega_V(\alpha).$$
We already know from Remark \ref{r4.1} that $L$ is
an additive group of ${\frak h}.$ By  Corollary \ref{ca} it is enough
to prove that $L$ is an even positive-definite lattice of rank $l.$ 

Note that the eigenvalue of $L(0)=L'(0)$ on  $\Omega_V(\alpha)$
is $\<\alpha,\alpha\>/2.$ It is clear that
$0\leq \<\alpha,\alpha\>/2\in \Z$ for all $\alpha\in L.$ That is, 
$L$ is even. If $\<\alpha,\alpha\>=0$ for some $\alpha\in L$ 
then $\Omega(\alpha)$ has weight zero.
Since $V_0$ is one dimensional and spanned by $\1$ we
immediately see that $\alpha=0.$ Thus $L$ is an even, positive-definite
lattice.

If the rank of $L$ is less than $l$ then there exists a non-zero $u\in \h$ 
such that $\<u, L\>=0.$ Thus $o(u)=u_0$ has zero eigenvalue
only on $V.$ As a result we see that 
$$\tr_Vo(u)o(v)q^{L(0)-c/24}=0$$
for all $v\in V_1.$ The proof of Theorem \ref{t1} shows that this
is not possible. As a result we see that $L$ is a positive-definite 
even lattice of rank $l.$ \qed

\end{document}